\documentclass[envcountsame, final]{llncs}

\usepackage{cite}
\usepackage{amsmath,amssymb,bbm}
\usepackage{mathrsfs}
\usepackage{subfig}
\usepackage{graphicx}
\usepackage{multirow}
\graphicspath{{./Figures/}}
\usepackage{microtype}
\usepackage{todonotes}
\usepackage{siunitx}
\let\six=\si

\def\th{\theta}

\def\si{\sigma} 
 
\def\ph{\varphi}

\def\SE{S\hspace{-0.08em}E}
\def\SO{SO}

\def\x{\times}

\let\on=\operatorname

\def\R{\mathbb{R}}

\def\Imm{\on{Imm}}

\newcommand{\ud}{\,\mathrm{d}}

\begin{document}

\title{Varifold-based matching of curves via Sobolev-type Riemannian metrics}
\author{Martin Bauer\inst{1}, Martins Bruveris\inst{2}, Nicolas Charon\inst{3}, Jakob M\o ller-Andersen\inst{1} \thanks{All authors contributed equally to the article.}}
\institute{Department of  Mathematics, Florida State University
\and Department of Mathematics, Brunel University London
\and CIS, Johns Hopkins University}

\maketitle

\begin{abstract}
Second order Sobolev metrics are a useful tool in the shape analysis of curves. In this paper we combine these metrics with varifold-based inexact matching to explore a new strategy of computing geodesics between unparametrized curves. We describe the numerical method used for solving the inexact matching problem, apply it to study the shape of mosquito wings and compare our method to curve matching in the LDDMM framework.
\end{abstract}

\section{Introduction}

Closed, unparametrized plane curves are used to represent the outline or shape of objects and as such they arise naturally in shape analysis and its applications~\cite{Srivastava2016}; these include medical imaging, computer animation, geometric morphometry and other fields. Analysis of shapes and their differences relies on the notion of distance between shapes. To define such a distance, we start from a Riemannian metric on the space of curves and compute its induced geodesic distance.

We consider in particular second order Sobolev metrics with constant coefficients. These are Riemannian metrics on $\on{Imm}(S^1,\R^d)$, the space of regular, parametrized curves and they are invariant under the reparametrization group $\on{Diff}(S^1)$. Hence they induce a Riemannian metric on the quotient space $\pi : \on{Imm}(S^1,\R^d) \to \on{Imm}(S^1,\R^d) /  \on{Diff}(S^1) \doteq B_{i}(S^1,\R^d)$ of unparametrized curves, whose elements are the shapes of objects one is interested in. To compute the geodesic distance between two shapes $\pi(c_0)$, $\pi(c_1)$ it is necessary to find minizing geodesics between the orbits $\pi(c_i) = c_i \circ \on{Diff}(S^1)$.

In previous work~\cite{BBHM2017} we achieved this by discretizing the diffeomorphism group $\on{Diff}(S^1)$ and its action on curves to obtain a numerical representation of the orbit $c_1 \circ \on{Diff}(S^1)$. Here we adopt a different approach. The varifold distance $d^{\on{Var}}$ between curves, defined in Sect.~\ref{sec:VarifoldDistance}, has reparametrizations in its kernel, meaning $d^{\on{Var}}(c_0, c_1) = d^{\on{Var}}(c_0, c_1 \circ \ph)$ and hence we can check equivalence of shapes or unparametrized curves via
\[
\pi(c_0) = \pi(c_1)
\quad\Leftrightarrow\quad
c_0 \in c_1 \circ \on{Diff}(S^1)
\quad\Leftrightarrow\quad
d^{\on{Var}}(c_0, c_1) = 0\,.
\]
Because the constraint $d^{\on{Var}}(c_0, c_1) = 0$ is difficult to encode numerically we relax it and solve an inexact matching problem instead. Given two curves $c_0$, $c_1$, to find the geodesic between the shapes they represent, we minimize
\[
E(c) + \lambda d^{\on{Var}}(c(1), c_1)\,,
\]
over all paths $c=c(t,\th)$ with $c(0) = c_0$, where $E$ is the Riemannian energy and $\lambda$ a coupling constant; see Sect.~\ref{sec:InexactMatching} for details.

We describe the numerical implementation in Sect.~\ref{implementation} and then apply in Sect.~\ref{experiments} the proposed framework to analyse the shape of mosquito wings~\cite{rohlf1984comparison} and to classify fish outlines~\cite{Mokhtarian1996} using geodesic distances and spectral clustering. We also compare the behaviour of Sobolev metrics with curve matching in the LDDMM framework.

\section{Mathematical background}

 \subsection{Sobolev metrics on shape space of curves}
 \label{sec:SobolevMetrics}
The space of smooth, regular curves with values in $\R^d$ is denoted by
\begin{align}
\Imm(S^1,\R^d)=\left\{c\in C^{\infty}(S^1,\R^d)\colon \forall \th \in S^1, c'(\th) \neq 0 \right\}\,,
\end{align}
where $\Imm$ stands for \emph{immersions}. We call such curves parametrized, because as maps from the circle they carry with them a parametrization; we will define the space $B_{i,f}(S^1,\R^d)$ of unparametrized curves in~\eqref{eq:Bif_def}. The space $\Imm(S^1,\R^d)$ is an open subset of the Fr\'echet space $C^\infty(S^1,\R^d)$ and therefore can be considered as a Fr\'echet manifold. Its tangent space $T_c\Imm(S^1,\R^d)$ at any curve $c$ is the vector space $C^\infty(S^1,\R^d)$ itself. 

We denote the Euclidean inner product on $\R^d$ by $\langle\cdot,\cdot\rangle$. Differentiation with respect to the curve parameter $\th \in S^1$ is written as $c_\theta=\partial_\theta c=c'$. For any fixed curve $c$, we denote differentiation and integration with respect to arc length by $D_s=\frac 1{|c_\theta|}\partial_{\theta}$ and $\mathrm ds=|c_\theta|\ud \theta$, respectively. A path of curves is a smooth map $c:[0,1]\to\Imm(S^1,\R^d)$ and we denote
the space of all paths by $\mathcal P = C^\infty([0,1], \Imm(S^1,\R^d))$. The velocity of a path $c$ is denoted by $c_t=\partial_t c=\dot c$. 
\begin{definition}
A second order Sobolev metric with constant coefficients is a Riemannian metric on the space $\Imm(S^1,\R^d)$ of parametrized curves of the form
\begin{equation}\label{eq:sobolev_metric}
G_c(h,k) = \int_{S^1} a_0\langle h,k \rangle+a_1\langle D_s h,D_s k \rangle+a_2\langle D_s^2 h,D_s^2 k \rangle \ud s \,,
\end{equation}
where $h,k \in T_c\Imm(S^1,\R^d)$ are tangent vectors, $a_j \in \R_{\geq 0}$ are constants and $a_0, a_2 > 0$. If $a_2=0$ and $a_1 > 0$, then $G$ is a first order metric and if $a_1=a_2=0$ it is a zero order or $L^2$-metric.
\end{definition}

Note that the symbols $D_s$ and $\ud s$ hide the nonlinear dependency of $G_c$ on the base point $c$. We use arc length operations in the definition of $G$ to ensure that the resulting metric is invariant under the action of the diffeomorphism group of $S^1$. The invariance property in turn allows us to define, using $G$, a Riemannian metric on the shape space of unparametrized curves.

The Riemannian length of a path $c\colon[0,1]\to\Imm(S^1,\R^d)$ is defined via
\begin{align}
L(c) = \int_0^1 \sqrt{G_{c(t)}(c_t(t),c_t(t))} \ud t\,.
\end{align}
The induced geodesic distance between two given curves $c_0$, $c_1$ of the Riemannian metric $G$ is then the infimum of the lengths of all paths connecting these two curves, i.e., 
\[
\on{dist}(c_0, c_1) = \inf \left\{L(c)\,:\, c \in \mathcal P,\, c(0)=c_0,\, c(1)=c_1\right\}\,.
\]
We can find critical points of the length functional by looking for critical points of the energy
\begin{equation} \label{eq:EnergyFunctional}
\begin{aligned}
E(c) &= \int_0^1 G_{c(t)}(c_t(t), c_t(t)) \ud t \,.
\end{aligned}
\end{equation}
The geodesic equation which corresponds to the first order condition for critical points, $DE(c) = 0$, is in the case of Sobolev metrics a partial differential equation for the function $c = c(t,\theta)$. Solutions of this equation are called geodesics, and are locally distance-minimizing paths. Local and global existence results for geodesics of Sobolev metrics were proven recently in \cite{Bruveris2014, Bruveris2014b_preprint, Michor2007} and they serve as the theoretical foundation of the proposed numerical framework: they
 tells us that the geodesic distance between two curves can always be realized by a  path between them, i.e. we can compute the geodesic distance by finding the energy-minimizing path.
 
\subsubsection*{Unparametrized curves.} Two curves that differ only by their parametrization represent the same geometric object. In the context of shape analysis it is therefore natural to consider them equal, i.e., we identify the curves $c$ and $c \circ \ph$, where $\ph$ is a reparametrization. We use as the reparametrization group the group,
\begin{equation*}
\on{Diff}(S^1) = \left\{ \ph \in C^\infty(S^1,S^1) \,:\, \ph' > 0 \right\}\,,
\end{equation*}
of smooth diffeomorphisms of the circle. This is an infinite-dimensional regular Fr\'echet Lie group \cite{Michor1997}. Reparametrizations act on curves by composition from the right, i.e., $c\circ \ph$ is a reparametrization of $c$. 

To define the quotient space of unparametrized curves we need to restrict ourselves to \emph{free immersions}, i.e. those upon which the diffeomorphism group acts freely. In other words
\[
c \in \Imm_f(S^1,\R^d)
\quad\Leftrightarrow\quad
\big(c \circ \ph = c \;\Rightarrow\; \ph = \on{Id}_{S^1}\big)\,.
\]
This restriction is necessary for technical reasons; in applications almost all curves are freely immersed. The space 
\begin{equation}
\label{eq:Bif_def}
 B_{i,f}(S^1,\R^d)=\Imm_f(S^1,\R^d) / \on{Diff}(S^1)\,,
\end{equation}
of unparametrized curves is the orbit space of the group action restricted to free immersions. This space is again a Fr\'echet manifold although constructing charts is nontrivial in this case \cite{Michor1991}. 

A Riemannian metric $G$ on $\on{Imm}(S^1,\R^d)$ is said to be \emph{invariant} with respect to reparametrizations if is satisfies
\[
G_{c \circ \varphi}( h \circ \varphi, k\circ \varphi) = G_c(h,k)\,,
\]
for all $\ph \in \on{Diff}(S^1)$. Sobolev metrics with constant coefficients are invariant with respect to reparametrizations and we have the following result concerning induced metrics on the quotient space.
\begin{theorem}\label{thm:shape_space} 
An Sobolev metric with constant coefficients on $\Imm(S^1,\R^d)$ induces a metric on $B_{i,f}(S^1,\R^d)$ such that the projection $\pi : \Imm_f(S^1,\R^d) \to B_{i,f}(S^1,\R^d)$ is a Riemannian submersion. 
\end{theorem}
The geodesic distance of the induced Riemannian metric on $B_{i,f}(S^1,\R^d)$ can be calculated using paths in $\Imm(S^1,\R^d)$ connecting $c_0$ to the orbit $c_1 \circ \on{Diff}(S^1)$, i.e., for $\pi(c_0),\pi(c_1) \in B_{i,f}(S^1,\R^d)$ we have,
\begin{equation*}
\on{dist}\big(\pi(c_0), \pi(c_1)\big) = \inf \left\{ L(c) \,:\, c \in \mathcal P,\, c(0) = c_0,\, c(1) \in c_1 \circ \on{Diff}(S^1)\right\} \,.	
\end{equation*}

In the same way the action of other groups can be factored out if the metric has a corresponding invariance property. 
We will consider later the action of the group of orientation-preserving Euclidean motions $\SE(d) = SO(d) \ltimes \R^d$,
\[
(A,b).c = A.c + b, \quad (A,b) \in \SO(d) \ltimes \R^d\,.
\]
The corresponding geodesic distance on the space $B_{i}(S^1,\R^d)/\SE(d)$
is then
\begin{equation*}
\on{dist}\big(\pi(c_0), \pi(c_1)\big) \!=\! \inf \left\{ L(c): c \in \mathcal P,\, c(0) = c_0,\, c(1) \in c_1 \circ \SE(d) \times \on{Diff}(S^1) \right\} \,.	
\end{equation*}

\subsection{Varifold distance on the space of curves}
\label{sec:VarifoldDistance}
A second construction of a distance on the space of curves arises from the framework of geometric measure theory by interpreting curves as currents or varifolds. These concepts go back to the works of Federer but have been recently revisited within the field of shape analysis as practical fidelity terms for diffeomorphic registration methods \cite{Glaunes2008,Charon_thesis}. Instead of relying on a quotient space representation, the core idea is to embed curves in a space of distributions. While the most general approach is explained in \cite{Charon2017}, the following gives a condensed overview adapted to the case of interest to this paper. 

Let $C_0(\R^d \times S^1)$ be the space of continuous functions vanishing at infinity. 
\begin{definition}
A varifold is an element of the distribution space $C_0(\R^d \times S^1)^*$. The varifold application $\mu: \ c \mapsto \mu_c$ associates to any curve $c \in \on{Imm}(S^1,\R^d)$ the varifold $\mu_c$ defined, for any $\omega \in C_0(\R^d \times S^1)$, by
\begin{equation}
 \label{def_varifold}
 \mu_c(\omega) = \int_{S^1} \omega\left(c(\theta),\frac{c'(\theta)}{|c'(\theta)|}\right) \!\ud s\,.
\end{equation}
\end{definition}
The essential property is that $\mu_c$ is actually independent of the parametrization in the sense that for any reparametrization $\ph \in \text{Diff}(S^1)$, one has $\mu_{c \circ \ph} = \mu_c$. Thus the map $c \mapsto \mu_c$ projects to a well-defined map from $B_{i,f}(S^1,\R^d)$ into the space of varifolds. 
Note however that the space of varifolds contains many other objects as well. 

Distances between curves can be defined by restricting a distance or pseudo-distance defined on the space of varifolds. A simple approach leading to closed form expressions is to introduce a \emph{reproducing kernel Hilbert space} (RKHS) of test functions and to use the corresponding kernel metric. Specifically, we consider kernels on $\R^d \times S^1$ of the form $k(x,u,y,v) \doteq \rho(|x-y|^2).\gamma(u \cdot v)$, i.e., $k$ is the product of a positive, continuous radial basis function $\rho$ on $\R^d$ and a positive, continuous zonal function $\gamma$ on $S^1$. To any such $k$ corresponds a RKHS $\mathcal H$ of functions, embedded in $C_0(\R^d \times S^1)$ with a dual metric $\langle \cdot,\cdot \rangle_{\on{Var}}$ on the corrsponding dual space $\mathcal H^\ast$ of varifolds. The reproducing kernel property implies---cf. \cite{Charon2017} for details---that for any curves $c_1,c_2$ we have
\begin{equation}
 \label{eq:metric_W_curves}
 \langle \mu_{c_1} , \mu_{c_2} \rangle_{\on{Var}} = \iint_{S^1\x S^1} \rho(|c_1(\theta_1) - c_2(\theta_2)|^2) \gamma\left(\frac{c_1'(\theta_1)}{|c_1'(\theta_1)|} \cdot \frac{c_2'(\theta_2)}{|c_2'(\theta_2)|} \right) \ud s_1 \ud s_2\,.
\end{equation}
Now, taking $d^{\on{Var}}(c_1,c_2) = \|\mu_{c_1} - \mu_{c_2}\|_{\on{Var}} = \langle \mu_{c_1} - \mu_{c_2}, \mu_{c_1} - \mu_{c_2} \rangle^{1/2}_{\on{Var}}$, the results of \cite{Charon2017} imply the following theorem.
\begin{theorem}
\label{thm:vardistance}
If $\rho$ and $\gamma$ are $C^1$ functions, $\rho$ is $c_0$-universal and $\gamma(1)>0$, then $d^{\on{Var}}$ defines a distance between any two closed, unparametrized, oriented and embedded curves. In addition, the distance is invariant with respect to the action of rigid transformations.   
\end{theorem}
Note that we require the stronger condition that the curves under consideration have to be embedded. On the bigger space of immersions, as considered in the previous section, the induced distance can degenerate. 

Invariance to rigid transformations means that for any $(A,b) \in \SE(d)$, $d^{\on{Var}}(A.c_1+b,A.c_2+b) = d^{\on{Var}}(c_1,c_2)$. It is also possible to construct distances that are invariant with respect to changes of orientation: to achieve this one simply selects kernels satisfying $\gamma(-t)=\gamma(t)$. 

Broadly speaking, the varifold distance \eqref{eq:metric_W_curves} results in a localized comparison between the relative positions of points and tangent lines of the the two curves, quantified by the choice of kernel functions $\rho$ and $\gamma$. As such, they do not derive from a Riemannian structre and there is no notion of geodesics in a shape space of curves. However, they provide a very efficient framework for defining and computing fidelity terms in relaxed matching problems as explained below.

 \subsection{Inexact matching on the shape space of curves}
 \label{sec:InexactMatching}
In this section we combine Sobolev metrics and varifold distances to compute geodesics on shape space via a relaxed optimization problem. Because reparametrizations lie in the kernel of the varifold distance,
\[
d^{\on{Var}}(c_0, c_1 \circ \ph) = \| \mu_{c_0} - \mu_{c_1 \circ \ph} \|_{\on{Var}}
= \| \mu_{c_0} - \mu_{c_1} \|_{\on{Var}} = d^{\on{Var}}(c_0, c_1)\,,
\]
we can reformulate the problem of finding geodesics as a constrained minimization problem,
\begin{equation}\label{Energy_Vari}
\on{dist}(\pi(c_0), \pi(c_1))^2 = 
\inf\left\{ E(c) \,:\, c \in \mathcal P,\, c(0) = c_0,\, d^{\on{Var}}(c(1),c_1)=0 \right\}\,.
\end{equation}
where $d^{\on{Var}}(c(1),c_1)$ is the varifold distance between the endpoint of the path $c(1)$ and the target curve $c_1$. 

Because it is difficult to numerically encode the constraint $d^{\on{Var}}(c(1),c_1) = 0$, we minimize the relaxed Lagrangian functional instead,
\begin{equation}\label{eq:Energy_Vari}
\on{dist}(\pi(c_0), \pi(c_1))^2 \approx \inf \left\{ E(c) +\lambda d^{\on{Var}}(c(1),c_1)^2
\,:\, c \in \mathcal P,\, c(0)=c_0\right\}\,.
\end{equation} 
If we minimize the relaxed functionals with an increasing sequence $\lambda \to \infty$, we would expect the minimizers to converge to the solution of the constrained minimization problem. For now we solve the problem with a fixed, large value of $\lambda$. This does not yield a geodesic with the correct endpoint, but the distance is small in practice.
In the future we plan to analyze the problem using an augmented Lagrangian approach to automatically choose a suitable value for $\lambda$.

Considering additionally the action of the Euclidean motion group $\SE(d)$
and using the invariance of the varifold distance under this group action we obtain the minimization problem
\begin{equation}\label{eq:Energy_Vari_Euc}
\inf\left\{ E(c) +\lambda d^{\on{Var}}(c(1),A.c_1+b)^2\,:\,
c \in \mathcal P,\, c(0)=c_0,\, (A,b) \in \SE(d) \right\}\,.
\end{equation}
to find geodesics on the space $B_{i.f}(S^1,\R^d)/\SE(d)$.

\section{Implementation}
\label{implementation}
\subsubsection*{The $H^2$-metric on spline-curves.}
\label{sec:H2Implementation}
In order to discretize the Riemannian energy term in the optimization problem \eqref{eq:Energy_Vari}, we discretize paths of curves using tensor product B-splines with $N_t \x N_\theta$ knots of orders $n_t = 2$ and $n_\theta = 3$,
\begin{equation}
\label{eq:TensorPath}
c(t, \th) = \sum_{i=1}^{N_t} \sum_{j=1}^{N_\th} c_{i,j} B_i(t) C_j(\th)\,.
\end{equation}
Here $B_i(t)$ are B-splines defined by an equidistant simple knot sequence on $[0,1]$ with full multiplicity at the boundary points, and $C_j(\theta)$ are defined by an equidistant simple knot sequence on $[0,2\pi]$ with periodic boundary conditions; for details see \cite{BBHM2017}. Note that the full multiplicity of boundary knots in $t$ implies
\begin{align*}
c(0, \th) &= \sum_{j=1}^{N_\th} c_{1,j} C_j(\th)\,, &
c(1, \th) &= \sum_{j=1}^{N_\th} c_{N_t,j} C_j(\th)\,.
\end{align*}
Thus the end curve $c(1)$ is given by the control points $c_{N_t,j}$. We approximate the integrals in the energy functional \eqref{eq:EnergyFunctional} using Gaussian quadrature with quadrature sites placed between knots.

Some notes on previous work: to solve the geodesic boundary value problem on shape space, we have proposed in \cite{BBHM2017} a method that involves discretizing the reparametrization group $\on{Diff}(S^1)$ using $B$-splines. However, the action of the reparametrization group is by composition, which does not preserve the $B$-spline space. To overcome this we added a projection step, where we project the composition $c \circ \ph$ back into the spline space. This has the disadvantage that the projection smoothes out details of the original curve, depending on how many control points are used and which parts of the curve are reparametrized. Furthermore, this methods requires a good choice of an initial path, which turned out to be a nontrivial obstacle in examples where the shapes under consideration are sufficiently different from each other. These considerations are our motivation to consider inexact matching with a varifold distance.

\subsubsection*{The varifold distance on spline curves.}
\label{sec:VarifoldDistanceImplementaion}
Our discretization of the varifold distance on spline curves builds on existing code for polygonal curves. Given two spline curves $c_{1}=\sum_{j=1}^{N_\th} c_{1,j} C_j$ and $c_{2}=\sum_{j=1}^{N_\th} c_{2,j} C_j$, a simple way of discretizing the varifold distance~\eqref{eq:metric_W_curves} is to approximate the splines by polygonal curves, i.e., choose sample vertices $v_{1,k}=c_{1}(\theta_k)$ and $v_{2,k}=c_2(\theta_k)$ with $\theta_k \in S^1$ for $1\leq k \leq P$ and $\theta_{P+1}=\theta_1$.
In future work we plan to calculate the varifold distance directly for spline curves, without the approximating step used here.

 Denoting the edge vectors $e_{1,k}=v_{1,k+1}-v_{1,k}$ and $e_{2,k}=v_{2,k+1}-v_{2,k}$, the inner product \eqref{eq:metric_W_curves} for the two polygonal curves $\tilde{c}_1$ and $\tilde{c}_2$ becomes
\begin{align*}
 \langle& \mu_{\tilde{c}_1} , \mu_{\tilde{c}_2} \rangle_{\on{Var}} = \\
 &\sum_{k,l=1}^{P} |e_{1,k}|.|e_{2,l}|.\gamma\left(\frac{e_{1,k}}{|e_{1,k}|} \cdot \frac{e_{2,l}}{|e_{2,l}|} \right) 
\iint_{[0,1]^2} \hspace{-1.5em} \rho(|v_{1,k}+t_1e_{1,k} - v_{2,l} -t_2 e_{2,l}|^2) \ud t_1 \ud t_2\,.
\end{align*}
In general, there is no closed form expressions for the double integral and hence we use a numerical approximation: we evaluate the integrand at the central point $(t_1,t_2) = (\frac 12, \frac 12)$, leading to the discrete approximation
\begin{align*}
\label{eq:metric_W_curves_discrete}
 \langle& \mu_{c_1} , \mu_{c_2} \rangle_{\on{Var}} \approx \\ &\sum_{k,l=1}^{P}
|e_{1,k}|.|e_{2,l}|.\gamma\left(\frac{e_{1,k}}{|e_{1,k}|} \cdot \frac{e_{2,l}}{|e_{2,l}|} \right)
.\rho \left(\left |\frac{v_{1,k}+v_{1,k+1}}{2} - \frac{v_{2,l}+v_{2,l+1}}{2} \right|^2 \right)
\end{align*}
and the corresponding expression for the distance $d^{\on{Var}}$. The total error resulting from both the polygonal approximation and the integral approximation can be shown to be of the order of $O(\max\{|\theta_{k+1}-\theta_{k}|\})$.

Additionally, we can compute the gradient of the discrete inner product with respect to, say, the spline coefficients $c_{1,j}$, usign the chain rule. Indeed, denoting $A$ the approximation of $\langle \mu_{c_1} , \mu_{c_2} \rangle_{\on{Var}}$,
\begin{equation*}
 \partial_{v_1} A = \partial_{x_1} A.\partial_{v_1} x_1 + \partial_{e_1} A.\partial_{v_1} e_1
\end{equation*}
where $x_{1,k} = (v_{1,k}+v_{1,k+1})/2$ is the edge midpoint. The gradients $\partial_{x_1} A$ and $\partial_{e_1} A$ are easily computed from the expression for $A$, while $\partial_{v_{1,l}}x_{1,k} = (\delta_{k-1}(l)+\delta_{k}(l))/2$ and $\partial_{v_{1,l}} e_{1,k} = \delta_{k-1}(l) - \delta_{k}(l)$. Finally, to obtain the gradient with respect the $c_{1,j}$ we apply the chain rule a second time, noting that $\partial_{c_{1,j}} v_{1,k} = C_j(\theta_k)$.

Our implementation includes many different choices of admissible kernel functions $\rho$ and $\gamma$, including the ones presented in \cite{Charon2017}. 

\subsubsection*{The inexact matching functional.}
With the discretization described above, the optimization problem \eqref{eq:Energy_Vari_Euc} becomes an unconstrained optimization problem for the control points $c_{i,j}$ the rotation matrix $A$ and the translation vector $b$. We choose a Limited-memory BFGS (L-BFGS) method to solve this problem, as implemented in the HANSO library~\cite{lewis2013nonsmooth} for Matlab, where we supply the formula for the gradient of the target function, see~\cite{BBHM2017} and~\cite{Charon2017} for the specific formulas. 

We initialize the optimization problem with the constant path $c(t,\theta) = c_0(\theta)$. Note that this overcomes one of the major drawbacks of the framework developed in \cite{BBHM2017}, which requires an initial path without singularities connecting the given curves $c_0$ and $c_1$.
To speed up the optimization we implemented a multigrid method, i.e., we first solve the geodesic problem with a coarser spline discretization and use the resulting optimal path to initialize the minimization of the original problem. A comparison of the resulting computation times can be seen in Table~\ref{table:computationtimes}. The obtained computation times are of the same order of magnitude as those of the SRV framework~\cite{Jermyn2011}\footnote{We used the publicly available Matlab implementation, which can be downloaded at \url{http://ssamg.stat.fsu.edu/software}.}.

\begin{table}
\begin{center}
\begin{tabular}{l c c c }
\hline\noalign{\smallskip}
{\bf Computation Times} & \phantom{i}L-BFGS\phantom{i} &\phantom{i}Multigrid\phantom{i} &\phantom{i}Aut. Diff.\phantom{i} \\
\noalign{\smallskip}\hline\noalign{\smallskip}
 Mosquito Wings & $1.9\six{\s}$ & $1.0\six{\s}$ & $3.0\six{\s}$  \\
 Surrey Fish & $3.1\six{\s}$ & $1.7\six{\s}$  & $4.7\six{\s}$     \\
 \hline
\end{tabular}
\end{center}
\caption{Average computation time (3.5 Ghz Core i7U) of the geodesic distance between mosquito wings (first line) and between shapes from the Surrey fish database (second line). The used methods are: L-BFGS without multigrid, with multigrid and 
with gradient calculation done via automatic differentiation (without multigrid). }
\label{table:computationtimes}
\end{table}
We also experimented with an implementation using automatic differentiation for the gradient calculation. While our implementation of the gradient is approximatively three times faster then the gradient computated with automatic differentiation, the resulting computation time for the optimization differed in average only by 51\% percent, see Table~\ref{table:computationtimes}. Automatic differentiation will allow us to implement, with little additional effort, a much wider class of Riemannian metrics, including curvature- and length-weighted metrics. See \cite{Michor2007,Bauer2014} for an overview of several Riemannian metrics on the space of curves. 

For our previous method~\cite{BBHM2017} we achieved a great speed-up of the optimization using a second order trust-region method---requiring us to compute the Hessian of the Riemannian energy. We tried this as well for this problem, but achieved no improvement in convergence. We speculate that this is due to the fact that the Hessian of the varifold distance is degenerate due to its kernel containing reparametrizations.

\section{Experiments}\label{experiments}
\subsubsection*{Influence of the weight $\lambda$.}
\begin{figure}
\centering
\includegraphics[width=\textwidth]{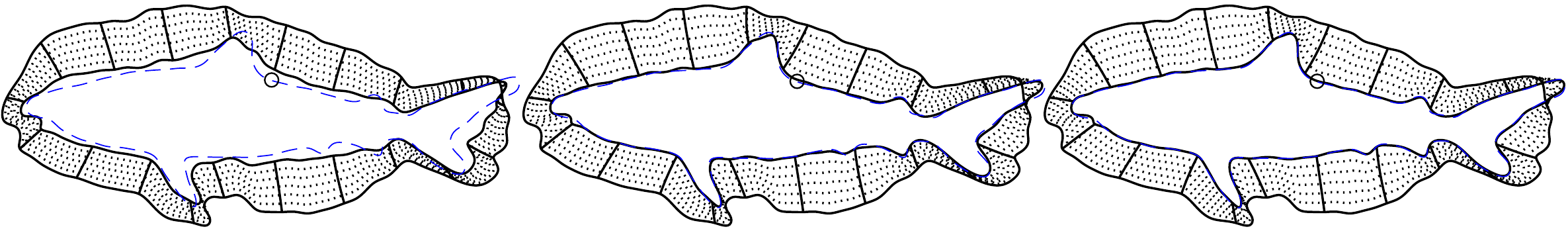}
\caption{Influence of the weight $\lambda$ on the quality of the matching: Minimizers of \eqref{eq:Energy_Vari_Euc} for $
\lambda = 0.3, 1$ and $5$. The target curve is depicted in blue.}
\label{fig:lambdainfluence}
\end{figure}

The weight $\lambda$ for the fidelity term  has a big influence on the quality of the matching. The solution of the optimization problem for different values of $\lambda$ is depicted in Fig.~\ref{fig:lambdainfluence}, and one can see that a good final matching requires a choice of a large enough $\lambda$. Choosing $\lambda$ too large, on the other hand, will make the functional too rigid, and the resulting deformation will be far from a geodesic. In the presence of noice, a large $\lambda$ might also result in overfitting. In future work we plan to investigate how to choose this parameter. 

\vspace{-1cm}
\subsubsection*{Towards a comparison with LDDMM curve matching.}
\begin{figure}
\centering
    \begin{tabular}{cccc}
\vspace{-.6cm}
\includegraphics[trim = 28mm 10mm 22mm 20mm ,clip,width=0.24\textwidth]{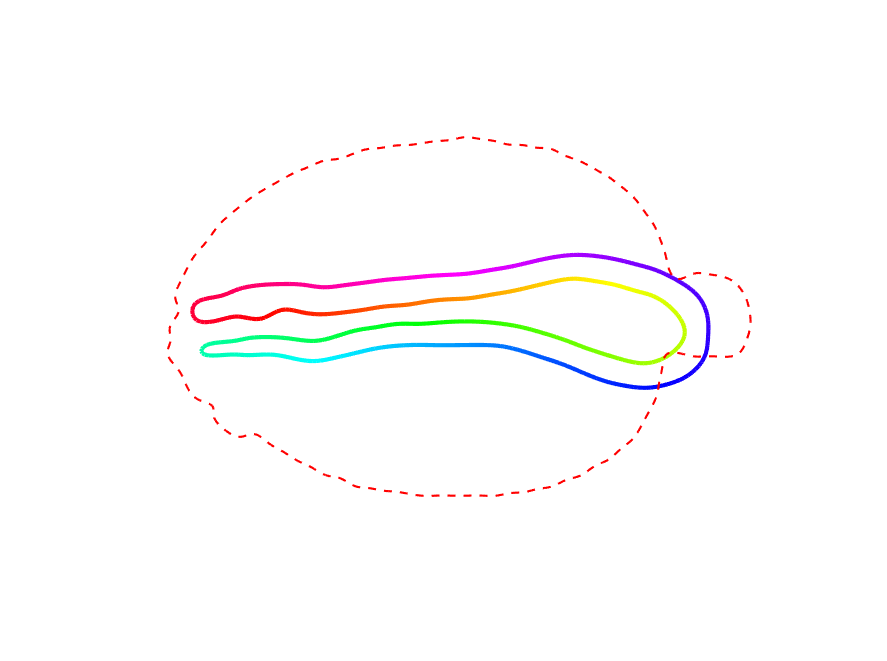} &
\includegraphics[trim = 28mm 11mm 22mm 20mm ,clip,width=0.24\textwidth]{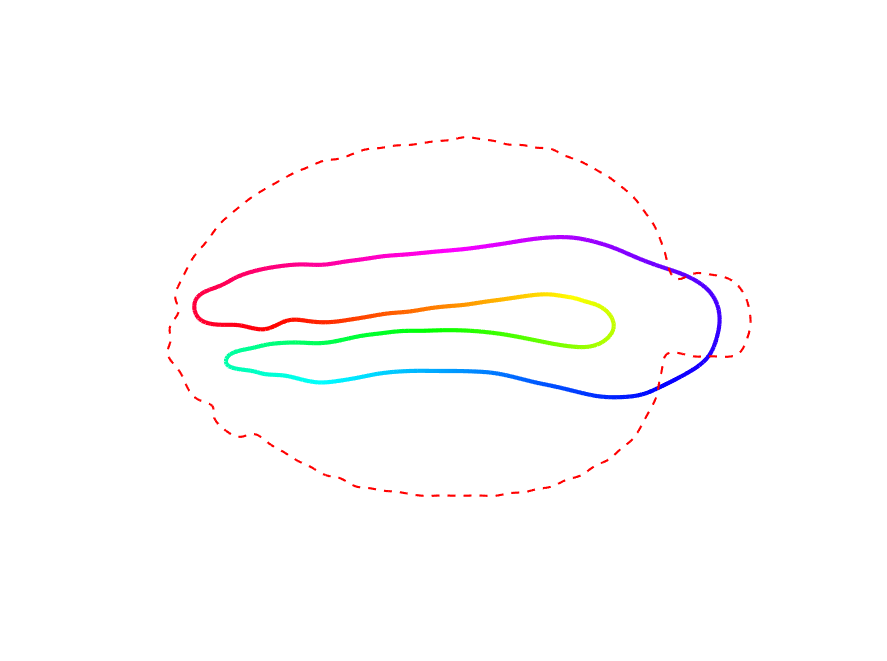} &
\includegraphics[trim = 28mm 11mm 22mm 20mm ,clip,width=0.24\textwidth]{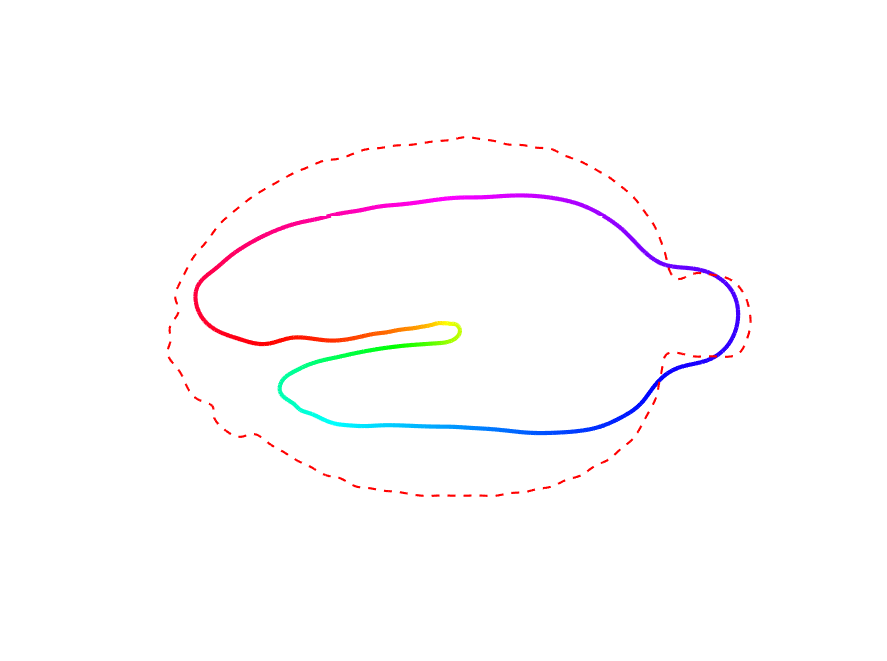} &
\includegraphics[trim = 28mm 11mm 22mm 20mm ,clip,width=0.24\textwidth]{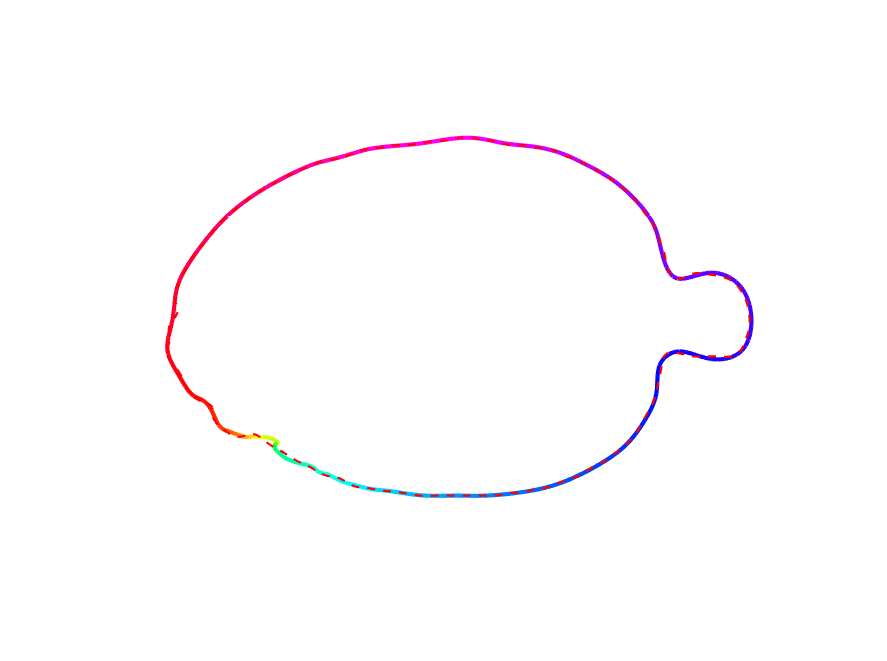} \\
\includegraphics[trim = 28mm 11mm 22mm 20mm ,clip,width=0.24\textwidth]{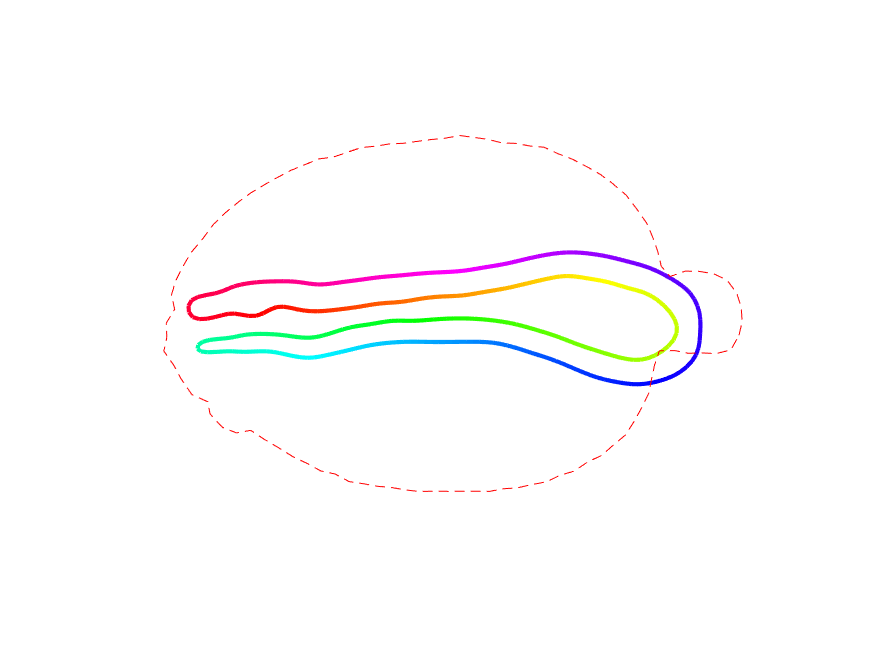} &
\includegraphics[trim = 28mm 11mm 22mm 20mm ,clip,width=0.24\textwidth]{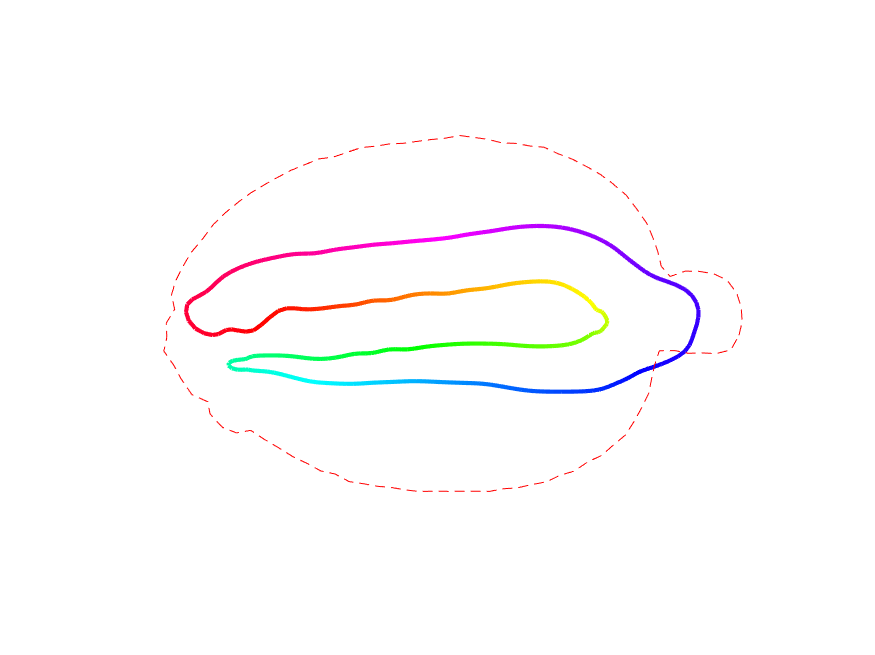} &
\includegraphics[trim = 28mm 11mm 22mm 20mm ,clip,width=0.24\textwidth]{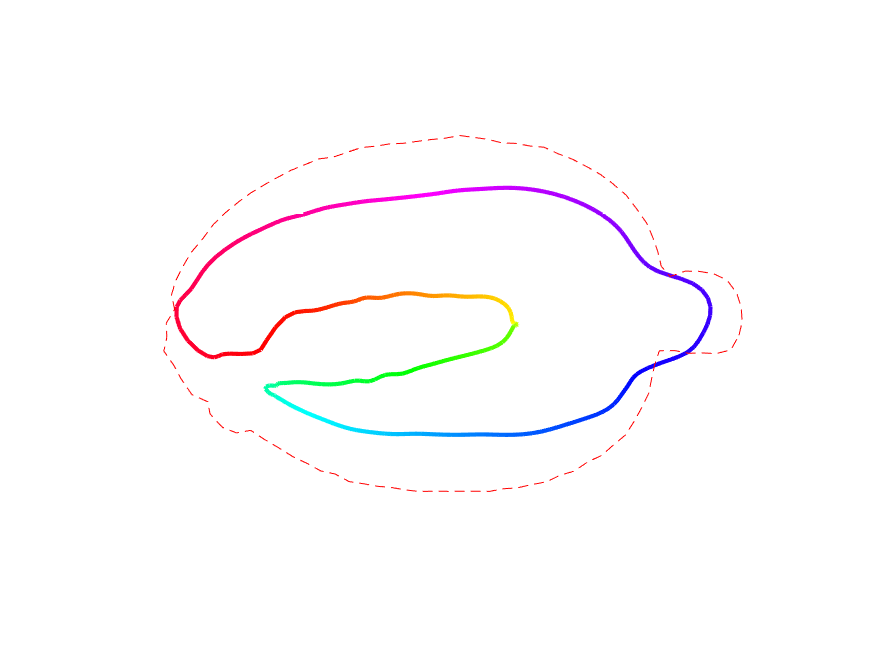} &
\includegraphics[trim = 27mm 11mm 22mm 20mm ,clip,width=0.24\textwidth]{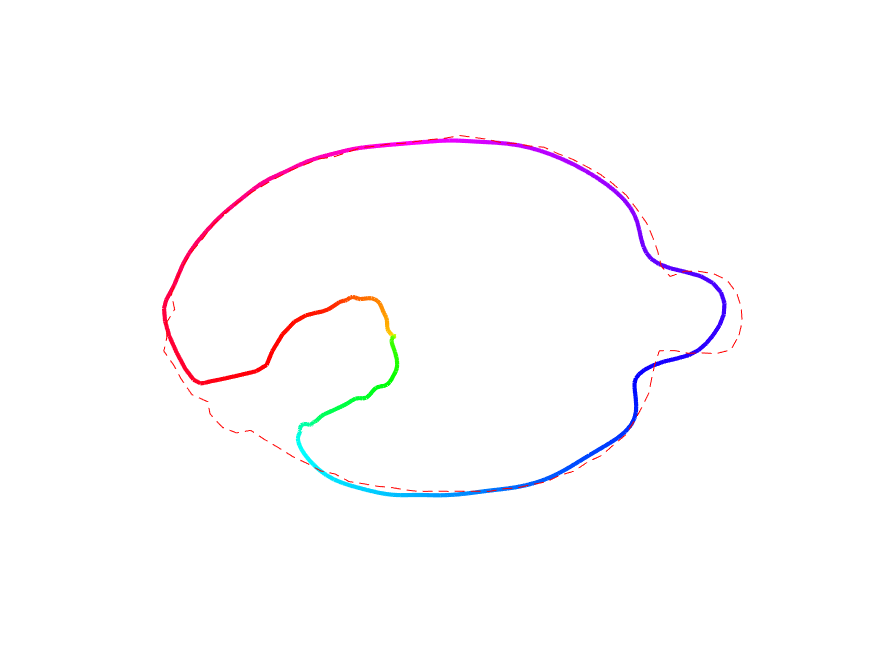} 
    \end{tabular}
    \caption{Geodesics at time steps $0,0.3,0.6,1$. First row: Inexact $H^2$-metric. Second row: LDDMM matching obtained using the algorithm of \cite{Charon2017}. Note that in the latter case, increasing the weight of the fidelity term or varying deformation scales do not in fact lead to significanty better matching results than shown here.} \label{fig:h2_vs_lddmm}
\end{figure}

Curve matching frameworks based on the LDDMM model such as \cite{Glaunes2008,Charon2017} are also formulated as relaxed optimizations involving the same varifold fidelity terms like \eqref{eq:Energy_Vari} with the difference that curve evolution is governed by an extrinsic and dense deformation of the plane and the metric on the shape space is now induced from a metric on the diffeomorphism group. Yet the parallels between the two formulations and algorithms should allow one to draw some insightful comparisons of the two models. Although this topic will need to be treated more extensively in future work, we show in Fig.~\ref{fig:h2_vs_lddmm} a simple example illustrating the difficulty, in the LDDMM setting, to generate a deformation that is able to stretch a thin structure in contrast with the intrinsic metric approach of this paper.     

\subsubsection*{Shape clustering.}
We next illustrate the discriminative power of our method for the problem of finding different clusters within a population of shapes. We focus on a small subset of $n=54$ shapes from the Surrey fish dataset and compute all the pairwise matchings between them. We then obtain a distance matrix given by the geodesic distance of the $H^2$-metric. Due to the asymmetry of inexact matching, we symmetrize the distance matrix a posteriori. In order to extract meaningful clusters, we use the spectral clustering framework presented in \cite{vonLuxburg2007}: the $p$-nearest neighbour graph is constructed based on the distance matrix (we use $p=12$ here) and the eigenvectors of the Jordan $\&$ Weiss normalized graph Laplacian are computed. Then each shape $i$ is mapped as the $i$-th row vector of the $n \times k$ matrix of the first $k$ eigenvectors and a $k$-means algorithm is used to separate those points into $k$ clusters.

\begin{figure}
\centering
    \begin{tabular}{cccc}
\includegraphics[trim = 15mm 10mm 15mm 0mm ,clip,width=0.24\textwidth]{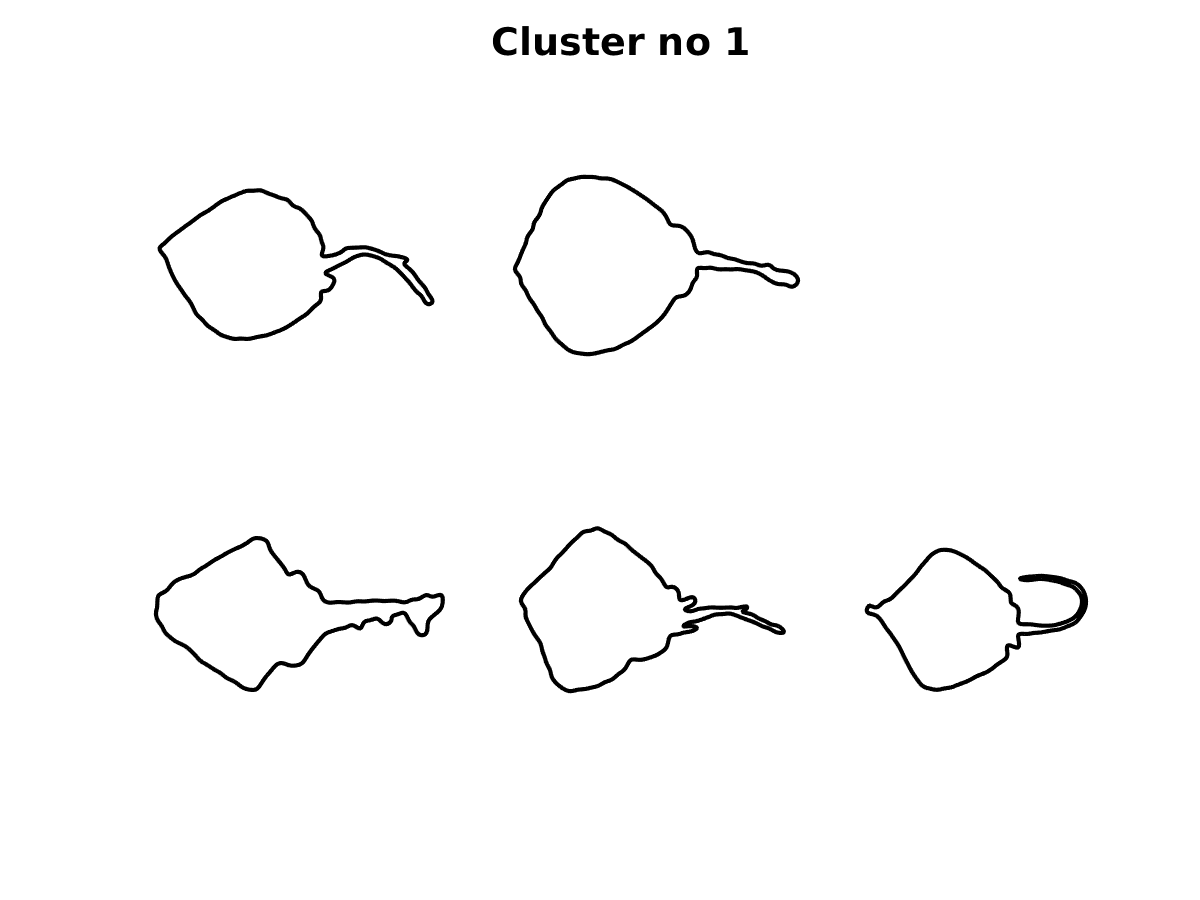} &
\includegraphics[trim = 15mm 10mm 15mm 0mm ,clip,width=0.24\textwidth]{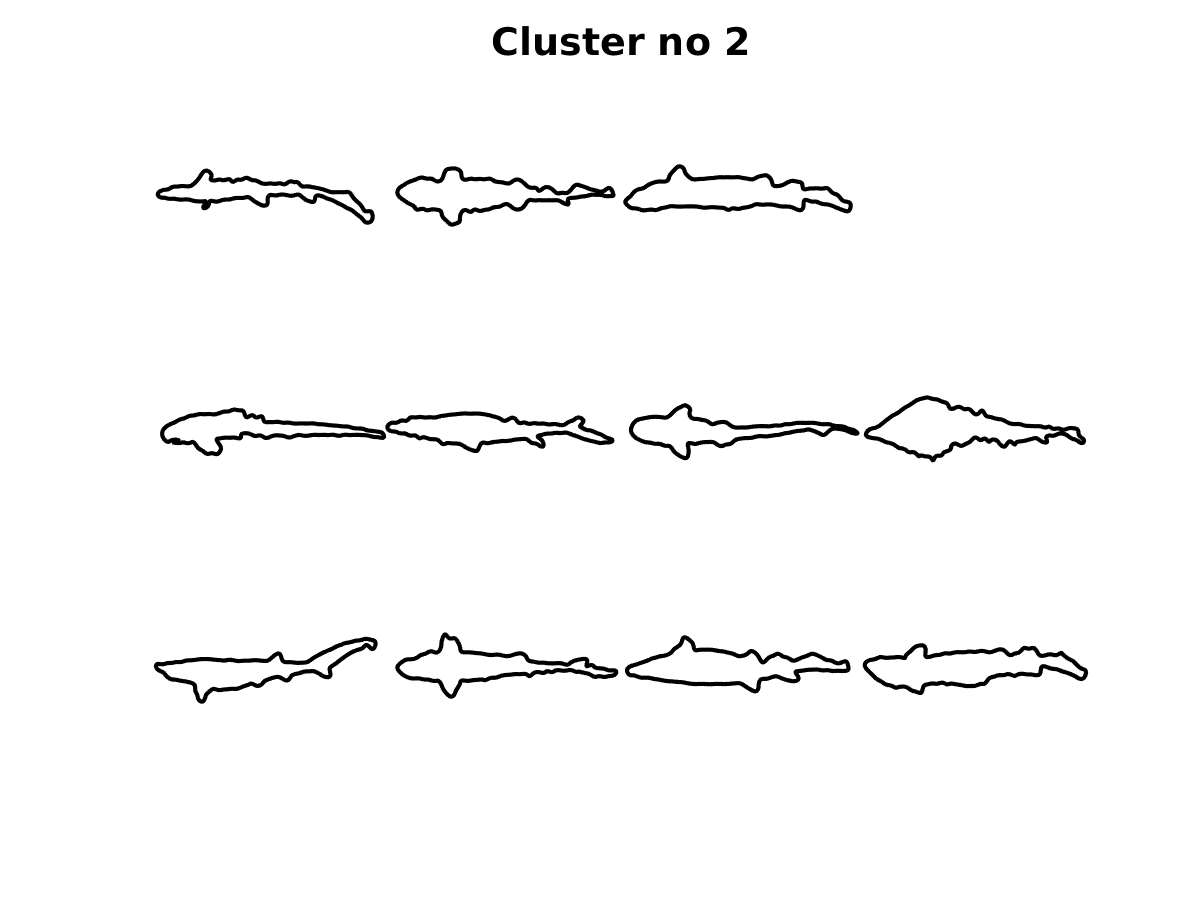} &
\includegraphics[trim = 15mm 10mm 15mm 0mm ,clip,width=0.24\textwidth]{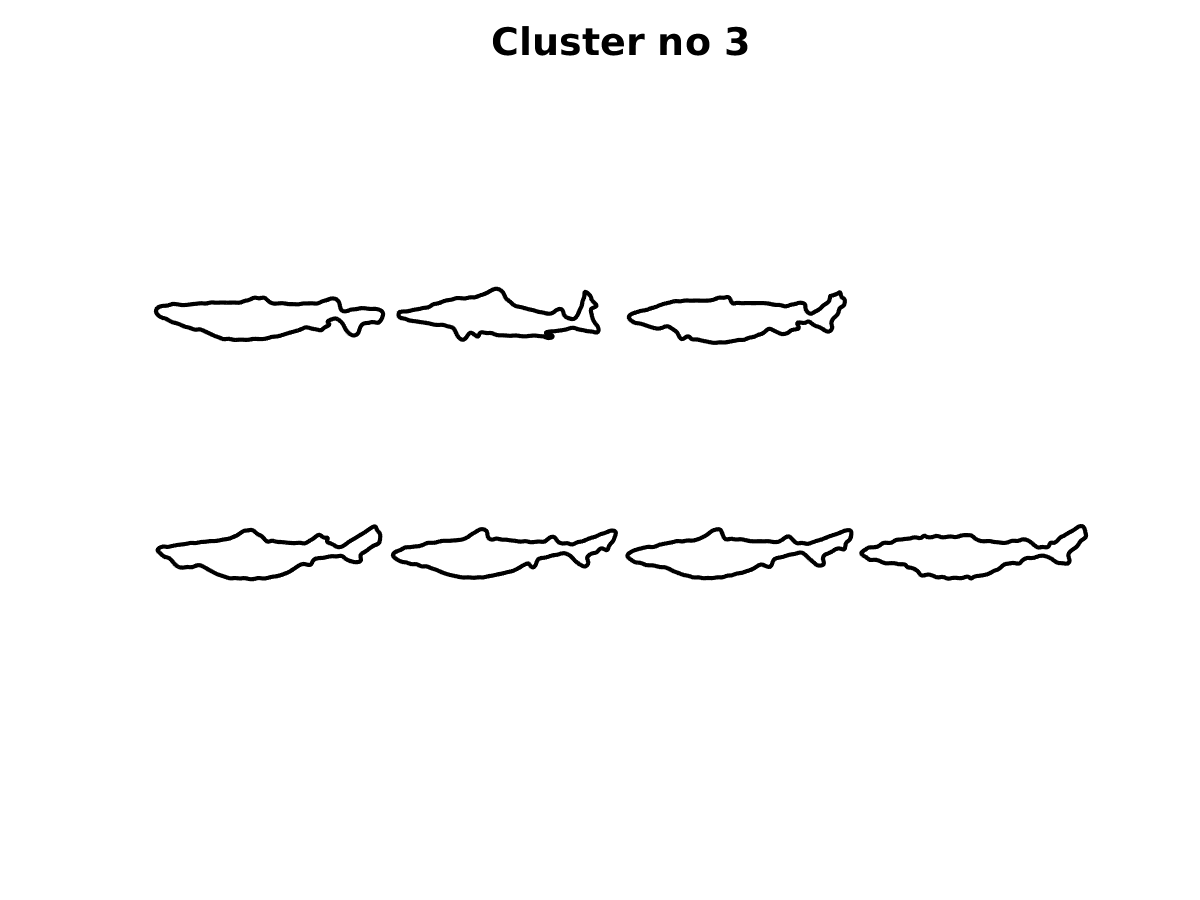} &
\includegraphics[trim = 15mm 10mm 15mm 0mm ,clip,width=0.24\textwidth]{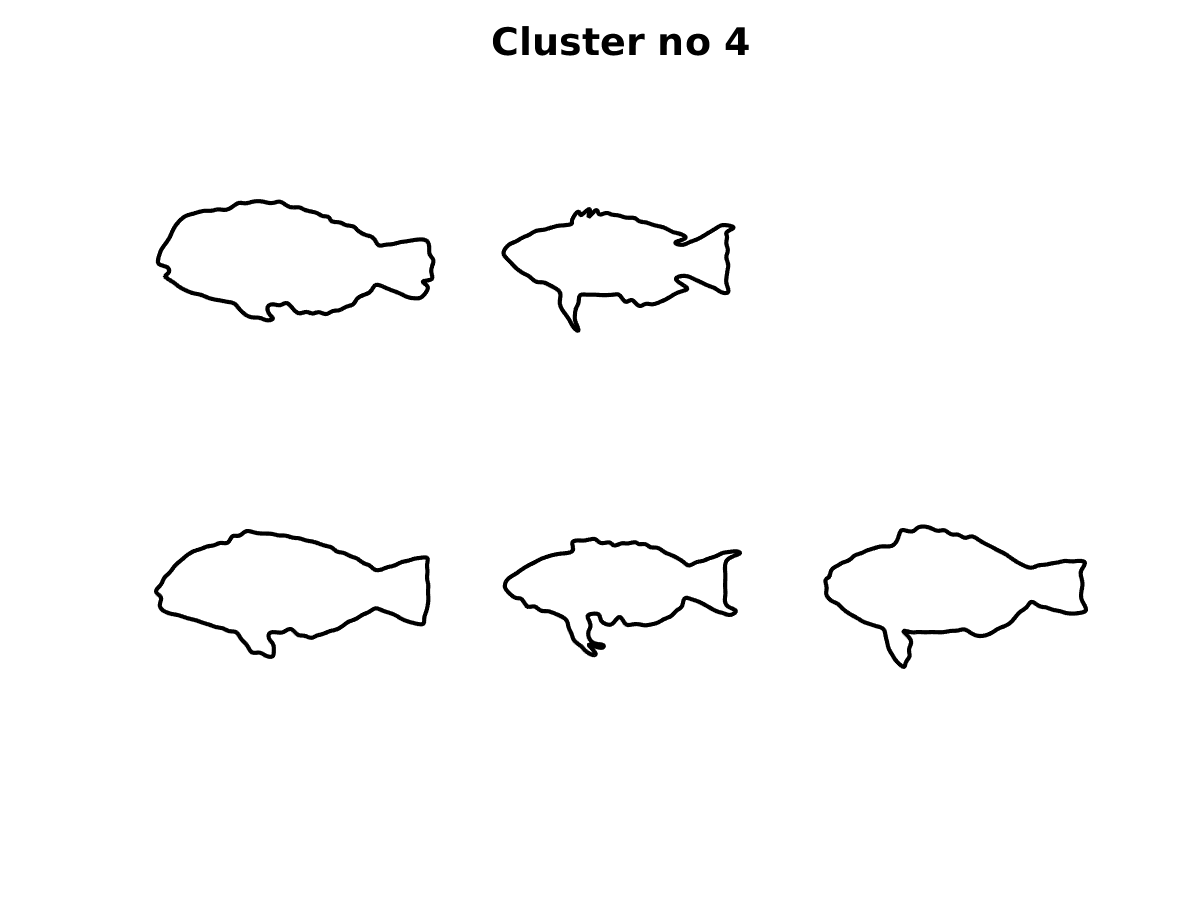} \\
\includegraphics[trim = 15mm 10mm 15mm 0mm ,clip,width=0.24\textwidth]{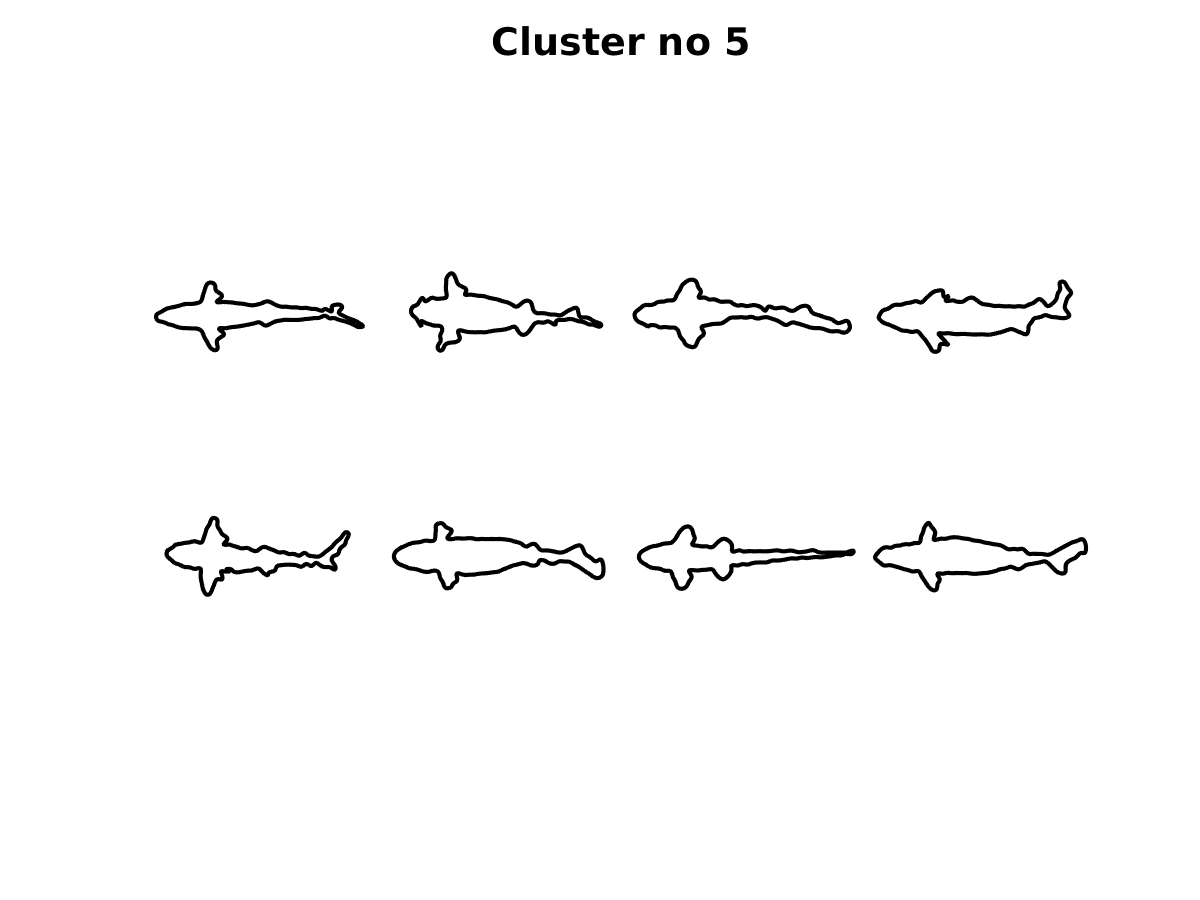} &
\includegraphics[trim = 15mm 10mm 15mm 0mm ,clip,width=0.24\textwidth]{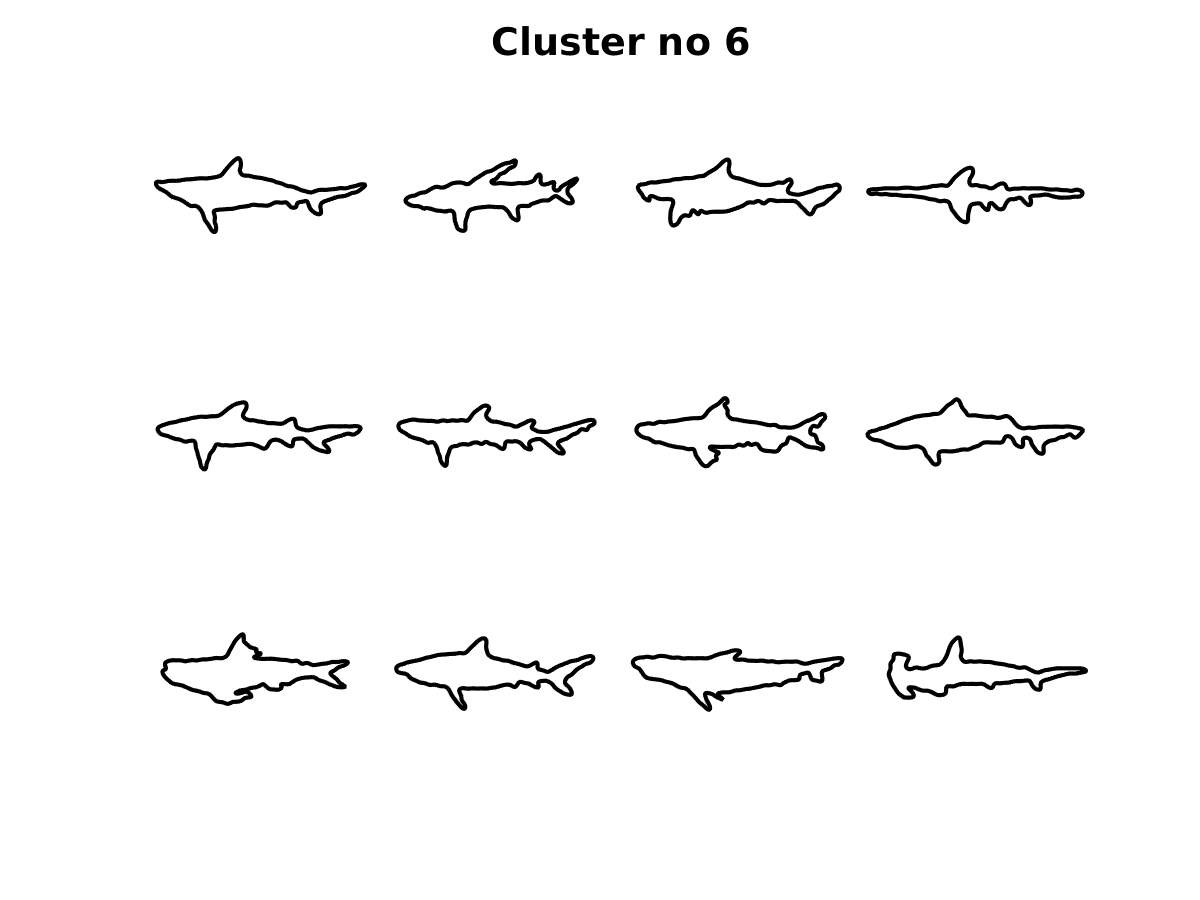} &
\includegraphics[trim = 15mm 10mm 15mm 0mm ,clip,width=0.24\textwidth]{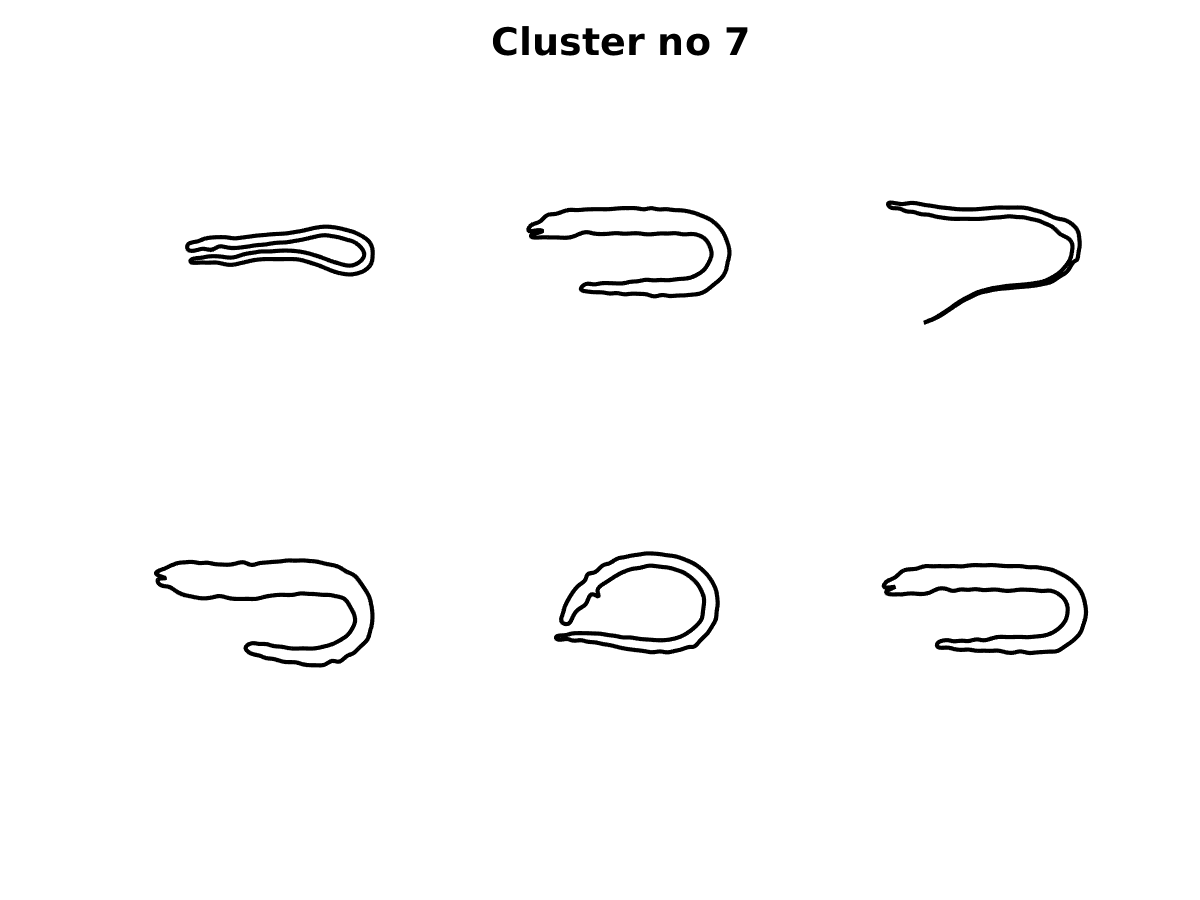} 
    \end{tabular}
    \caption{The results of the cluster analysis for 54 shapes from the Surrey fish dataset obtained from the spectral clustering method.} \label{fig:clustering}
\end{figure}

The results of this approach for $k=7$ clusters is shown in Fig.~\ref{fig:clustering}. Overall, up to a few exceptions, the method is able to discriminate well between the different classes of this particular population with an accuracy comparable to using the LDDMM framework for measuring shape distances but with significantly faster computation of the distance matrix.   

\subsubsection*{Mosquito wings.}
\begin{figure}
\centering
\begin{tabular}{cc}
\vspace{1cm}
\multirow{-2}{*}{}{\smash{\raisebox{-0.57\height}{\hspace{-0.2cm}
\includegraphics[width=0.53\textwidth]{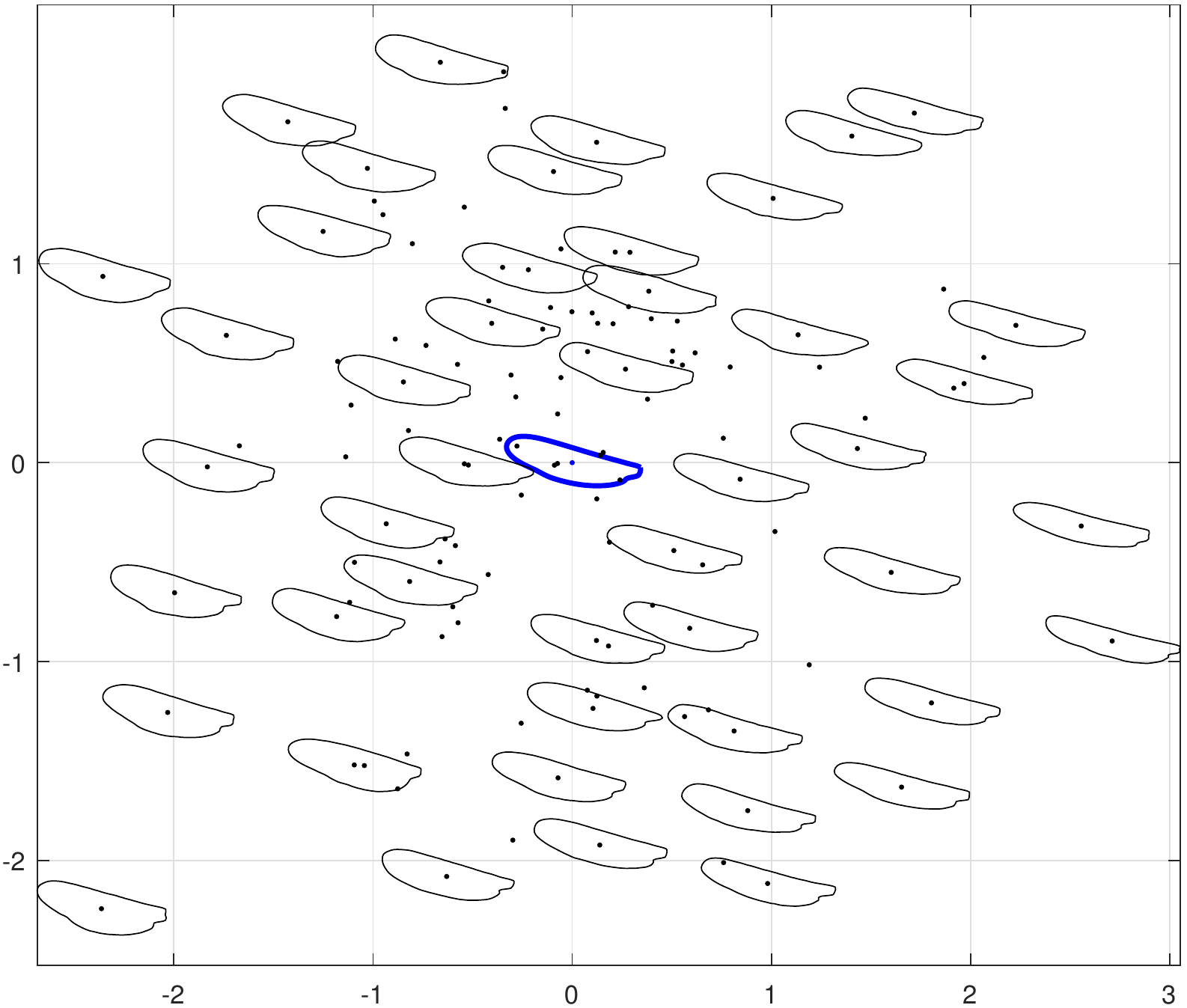}
\hspace{0.3cm}}}} 
&
\includegraphics[width=0.38\textwidth]{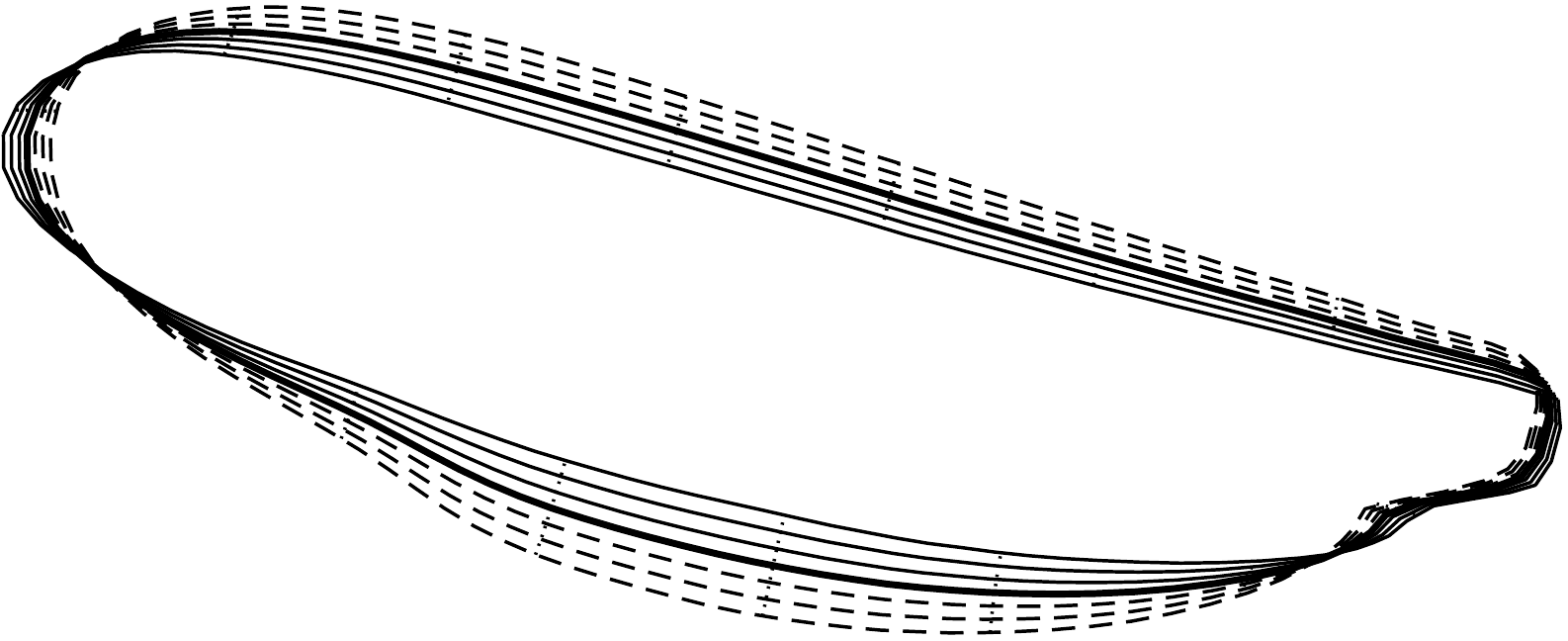} \\
& \includegraphics[width=0.38\textwidth]{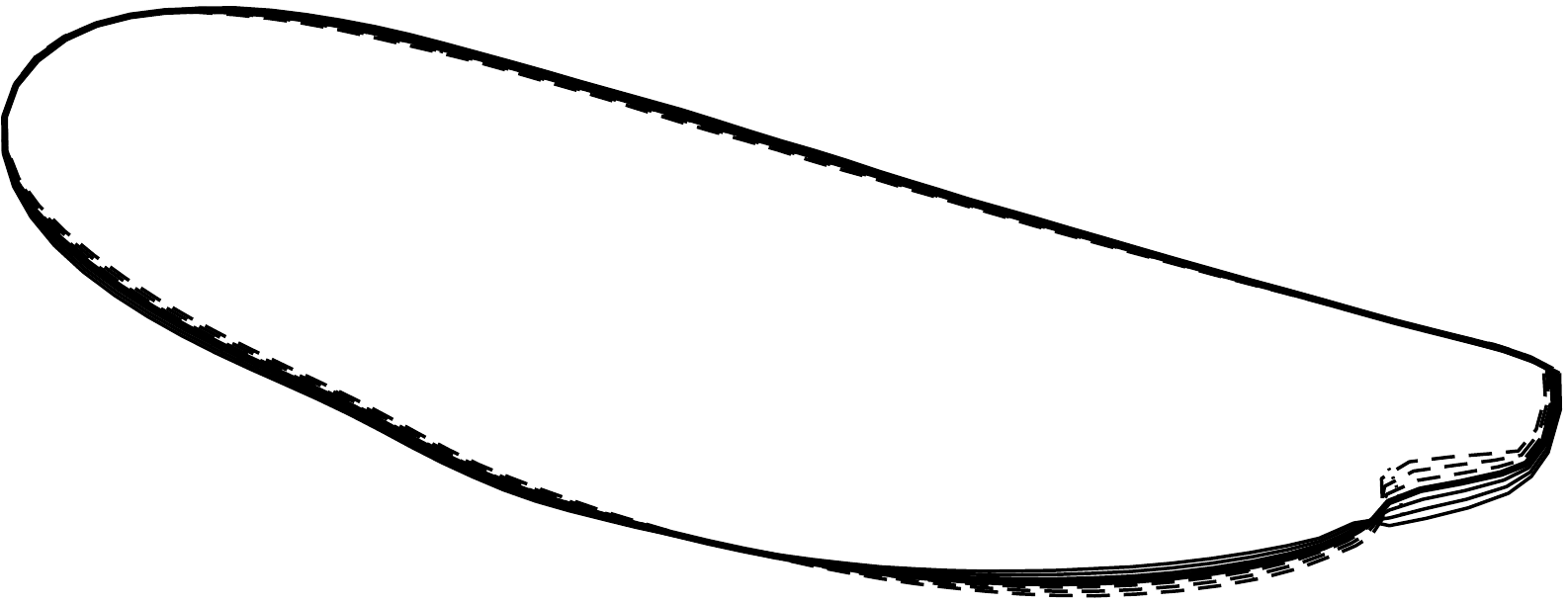}\vspace{.4cm}
\end{tabular}
\caption{Left: Wings projected to the plane in the tangent space of the mean, spanned by the first two principal directions. Right: Geodesics from the mean in the two first principal directions.}
\label{fig:wingsPCA}
\end{figure}

Finally we want to demonstrate our numerical framework by providing a simple analysis of a set of mosquito wings. 
The data consists of the boundary curves $c_1, \dots, c_{126}$ of 126 mosquito wings. The acquisition of the data is described in the article \cite{rohlf1984comparison}, where the authors analyzed the data using a polar coordinate system to describe each wing via the distance function from its centroid. Using inexact geodesic matching we use as a template the Karcher mean $\bar{c}$ of the data, with respect to the Sobolev metric. In the tangent space of the mean we represent each curve $c_j$ by the initial velocity $v_j = \on{Log}_{\bar{c}(c_j)}$ to the geodesic connecting $\bar{c}$ and $c_j$. We then do a PCA with respect to the inner product $G_{\bar{c}}$ of the data in the linear tangent space. Fig.~\ref{fig:wingsPCA} depicts the data set after projecting onto the subspace spanned by the first two principal components, and the geodesic in the direction of these to directions. It seems to suggest that the two principal directions control the thickness of the wings, and the depth of the fold at the end, respectively.

\section{Conclusions}
In this article we present a new numerical method to compute geodesics for second order Sobolev metrics. The proposed algorithm is based on previous work on Sobolev metrics \cite{BBHM2017} and the varifold distance \cite{Charon2017}.
Since reparametrizations are in the kernel of the varifold distance we avoid having to discretize the reparametrization group $\on{Diff}(S^1)$ to solve the geodesic boundary value problem on shape space. This allows us to  
overcome certain problems of the framework presented in  \cite{BBHM2017}. Furthermore this new approach is better suited for generalization to shape spaces of unparametrized surfaces. Additionally, since we are now using an L-BFGS method---this only requires the computation of the gradient but not of the Hessian---it will be possible to generalize this framework with little additional effort to a much wider class of metrics. We plan to follow these lines of research in future work and use it 
to investigate methods for a data-driven choice of a Riemannian metric for applications in shape analysis.

\bibliographystyle{splncs03}

\end{document}